\documentclass{ifacconf}

\usepackage{graphicx}      
\usepackage{natbib}        
\usepackage{amsmath, amssymb}
\newtheorem{theorem}{Theorem}
\newtheorem{lemma}{Lemma}

\newtheorem{assumption}{Assumption}

\newtheorem{remark}{Remark}
\usepackage{csquotes}
\def\R{\mathbb{R}}

\makeatletter
\DeclareRobustCommand{\qed}{%
	\ifmmode 
	\else \leavevmode\unskip\penalty9999 \hbox{}\nobreak\hfill
	\fi
	\quad\hbox{\qedsymbol}}
\newcommand{\openbox}{\leavevmode
	\hbox to.77778em{%
		\hfil\vrule
		\vbox to.675em{\hrule width.6em\vfil\hrule}%
		\vrule\hfil}}
\newcommand{\qedsymbol}{\openbox}
\newenvironment{proof}[1][\proofname]{\par
	\normalfont
	\topsep6\p@\@plus6\p@ \trivlist
	\item[\hskip\labelsep\itshape
	#1.]\ignorespaces
}{%
	\qed\endtrivlist
}
\newcommand{\proofname}{Proof}
\makeatother

\begin{document}
\begin{frontmatter}

\title{Predictor-Feedback Stabilization of Globally Lipschitz Nonlinear Systems with State/Input Quantization\thanksref{footnoteinfo}}

\thanks[footnoteinfo]{Funded by the European Union (ERC, C-NORA, 101088147). Views and 
opinions expressed are however those of the authors only and do not necessarily reflect those 
of 	the European Union or the European Research Council Executive Agency. Neither the 
European Union nor the granting authority can be held responsible for them.}

\author[1]{Florent Koudohode}
\author[1]{Nikolaos Bekiaris-Liberis} 
	
\address[1]{Department of Electrical and Computer Engineering, Technical 
		University of Crete, University Campus, Akrotiri, Chania, Greece}

\begin{abstract}               
We develop a switched nonlinear predictor-feedback control law to achieve global asymptotic stabilization for nonlinear systems with arbitrarily long input delay, under state quantization. The proposed design generalizes the nonlinear predictor-feedback framework by incorporating quantized measurements of both the plant and actuator states into the predictor state formulation. Due to the mismatch between the (inapplicable) exact predictor state and the predictor state constructed in the presence of state quantization, a global stabilization result is possible under a global Lipschitzness assumption on the vector field, as well as under the assumption of existence of a globally Lipschitz, nominal feedback law that achieves global exponential stability of the delay/quantization-free system. To address the constraints imposed by quantization, a dynamic switching strategy is constructed, adjusting the quantizer's tunable parameter in a piecewise constant manner—initially increasing the quantization range, to capture potentially large system states and subsequently refining the precision to reduce quantization error. The global asymptotic stability of the closed-loop system is established through solutions estimates derived using backstepping transformations, combined with small-gain and input-to-state stability arguments. We also extend our approach to the case of input quantization.
\end{abstract}


\end{frontmatter}

\section{Introduction}
The challenge of compensating long input delays in nonlinear control systems has been addressed through predictor-based feedback control design techniques. These methods ensure stability and robustness under ideal continuous-in-time implementations, as demonstrated in, for example, the works by \cite{krstic2009input,bekiaris2013nonlinear,karafyllis2017predictor}. However, real-world applications often involve digital implementation effects, such as sampling and quantization, which can degrade performance or even destabilize the system if not properly addressed, see, for example, \cite{mazenc2016predictor,karafyllis2017predictor}. Consequently, ensuring stability under such digital implementation constraints for nonlinear systems with long input delays is a challenging problem of theoretical and practical significance.

To address these challenges, previous studies investigate various digital effects in control systems under predictor feedback. The impact of sampling on measurements and control inputs is studied in \cite{karafyllis2011nonlinear,zhu,Selivanov2016, battilotti2019continuous, weston2018sequential}; whereas relevant, event-triggered predictor-based designs are developed in \cite{gonzalez2019event,mazenc2022event,nozari2020event}, and \cite{sun2022predictor}. Quantization effects in delay systems are explored for linear time-delay systems with saturation by \cite{fridman2009control} and for linear networked control systems by \cite{liu2015dynamic}, for logarithmic quantizers within an event-based control framework by \cite{garcia2012model}, for nonlinear systems with (small) measurement delay by \cite{liberzon2006quantization}, and for nonlinear systems with state delay by \cite{di2020practical,DIFERDINANDO2024111567}. Practical stabilization under static quantization is also investigated in \cite{espitia2017stabilization} and \cite{tanwani2016input} for first-order hyperbolic PDE systems, while \cite{selivanov2016predictor} and \cite{katz2022sampled} address parabolic PDE systems. Dynamic quantizers, which can adjust online their range and precision, achieving global asymptotic stabilization, are considered for first-order hyperbolic systems in \cite{bekiaris2020hybrid}, for linear systems with input delay in \cite{fkoudohode2024}, and for general infinite-dimensional discrete-time systems in \cite{wakaiki2024stabilization}. There exists no result addressing the problem of simultaneous compensation of long input delay and state/input quantization for nonlinear systems.

In this paper, we develop a novel, switched nonlinear predictor-feedback control law that achieves global asymptotic stabilization (in the supremum norm of the actuator state), for nonlinear systems with state quantization. The design approach that we introduce relies on two main ingredients—a quantized version of nonlinear predictor feedback, in which quantized measurements of both plant and actuator states are involved in the predictor state formula; and a dynamic switching strategy that properly adjusts the quantizer's tunable parameter in a piecewise constant manner (inspired by \cite{brockett2000quantized, liberzon2003hybrid} for the case of delay-free systems). In particular, the control design involves two main stages, referred to as \enquote{zooming out} and \enquote{zooming in} phases. In the zooming-out phase, in which the system operates in open loop, the quantization range expands until the system's infinite-dimensional state gets within the quantizer's range, while in the zooming-in (closed-loop) phase, the range is refined to reduce quantization error. Stability of the closed-loop system is established using backstepping transformations, in combination with derivation of solutions estimates and utilization of small-gain and input-to-state stability (ISS) arguments. Due to the mismatch between the (inapplicable) exact predictor state and the predictor state constructed in the presence of state quantization, a global stabilization result, which is established here, is possible under a global Lipschitzness assumption on the vector field, as well as under the assumption of existence of a globally Lipschitz, nominal feedback law that achieves global exponential stability of the delay/quantization-free system (see also, for example, \cite{bresch2014delay}). We also study existence and uniqueness of solutions under a locally Lipschitz assumption on the quantizers, which enables us to employ the results from \cite{karafyllis2021input}. While such an assumption on the quantizer may appear as restrictive, even though it is in fact employed only for analyzing existence/uniqueness and the stability estimates derived do not explicitly depend on it, it does not seem to be straightforward to remove it, at least in the case of state quantization in which the quantized, full actuator state is employed in an infinite-dimensional feedback law (see also \cite{bekiaris2020hybrid, fkoudohode2024}). We extend our approach to the case of input quantization as well. 

The paper is organized as follows. Section~\ref{probFormulation} introduces the system dynamics and the switched predictor-feedback design. Section~\ref{stabStateQuantization} establishes global asymptotic stability of the closed-loop system under state quantization, while Section~\ref{inputquantization} extends the stability results to input quantization. Section~\ref{conclusion} presents conclusions and future research directions.

{\em Notation:} We denote by \( L^{\infty}(A ; \Omega) \) the space of measurable and bounded functions defined on \( A \) and taking values in \( \Omega \). For a given \( D > 0 \) and a function \( u \in L^{\infty}([0, D]; \R) \), we define the \( L^\infty \)-norm of \( u \) as  
$\|u\|_{\infty} = \operatorname{ess \, \sup}_{x \in [0, D]} |u(x)|,$ where $\operatorname{ess \, \sup}$ is the essential supremum. For a real number \( h \in \mathbb{R} \), we define its integer part as  $\lfloor h \rfloor = \max \{ k \in \mathbb{Z} : k \leq h \}.$ The state space \( \mathbb{R}^n \times L^{\infty}([0, D]; \mathbb{R}) \) is equipped with the norm $\|(X,u)\| = |X| + \|u\|_{\infty},$ where \( X \in \mathbb{R}^n \) and \( u \in L^{\infty}([0, D]; \mathbb{R})\). We denote by \( AC(\mathbb{R}_+, \mathbb{R}^n) \) the set of all absolutely continuous functions \( X: \mathbb{R}_+ \to \mathbb{R}^n \). Let \( I \subseteq \mathbb{R} \) be an interval. The set of all piecewise right-continuous functions \( f: I \to J \) is denoted by \( \mathcal{C}_{\rm rpw}(I, J) \) (see also \cite{espitia2017stabilization,KKnonlocal}). 
\section{ Problem Formulation and Control Design}\label{probFormulation}
\subsection{ Nonlinear Systems With State Quantization} 
We consider the following nonlinear system
\begin{equation} \label{nonlinear_delay_system}
	\dot{X}(t)=f(X(t),U(t-D)),
\end{equation}
where $D>0$ is input delay, $t \geq 0$ is time variable, $X \in \Rset^{n}$ is state, $U$ is scalar control input, and $f:\R^n\times \R\to \R^n$ is vector field. System \eqref{nonlinear_delay_system} can be alternatively represented as follows 
\begin{align}
	\label{pde_representation}	\dot{X}(t)&=f(X(t),u(0,t)),\\
	\label{pde_representation01}	u_{t}(x, t)&=u_{x}(x,t), \\
	\label{pde_representation1}	u(D, t)&=U(t),
\end{align}
by setting $u(x, t)=U(t+x-D),$ where $ x\in [0, D]$ and $u$  is the transport PDE actuator state. 
We proceed from now on with representation \eqref{pde_representation}--\eqref{pde_representation1} as it turns out to be more  convenient for control design and analysis. With the backstepping transformations (direct and inverse),
\begin{align}
	w(x, t)&=u(x,t)-\kappa(p(x,t)),\\
	u(x,t)&=w(x,t)+\kappa(\pi(x,t)),\label{backstepping_direct_transformation}
\end{align}
in \cite{krstic2009input}, system \eqref{pde_representation}--\eqref{pde_representation1} is transformed into 
\begin{align}\label{targetsystem_without_quantizer}
	\dot{X}(t)&=f(X(t),\kappa(X(t))+w(0,t)), \\
	\label{targetsystem_without_quantizer2}	w_{t}(x, t)&=w_{x}(x, t), \\
	w(D, t)&=U(t)-U_{\rm nom}(t),\label{targetsystem_without_quantizer3}
\end{align} where $U_{\rm nom}(t)$ is the nominal predictor feedback defined as follows \begin{align}\label{nominalU}
	U_{\rm nom}(t)&=\kappa\left(P(t)\right),
\end{align} with $P(t)=p(D,t),$ where $p$ and $\pi$ are predictor variables, represented by the following integral equations
\begin{align}\label{pxt}
	p(x, t)&=\int_0^x f(p(\xi, t), u(\xi, t)) d \xi+X(t),\\
	\label{pixt}	\pi(x, t)&=\int_0^x f(\pi(\xi, t),\kappa(\pi(\xi, t))+w(\xi, t)) d \xi+X(t),
\end{align}
with $p,\pi:[0, D] \times \mathbb{R}_{+} \rightarrow \mathbb{R}^n$. 
We make the following assumptions. 
\begin{assumption}\label{globalLipschitz}
	The function $f: \R^n\times\R\to\R^n,$ which satisfies $f(0,0)=0,$ is continuously differentiable and globally Lipschitz, and thus, there exists $L>0$ such that $ \forall  u_1,u_2\in \R$ and $ \forall X_1,X_2\in \R^n,$
	\begin{equation}\label{lipschitzf}
		|f(X_1,u_1)-f(X_1,u_2)|\le L|X_1-X_2|+L|u_1-u_2|.
	\end{equation}
\end{assumption}

\begin{assumption}\label{issODE}
	The system $\dot{X} = f(X, \kappa(X))$ is globally exponentially stable, where $\kappa:\R^n\to \R$, satisfying $\kappa(0)=0,$ is a continuously differentiable, globally Lipschitz function, and hence, there exists a constant $\kappa_0 > 0$ such that for all $p, \pi \in \mathbb{R}^n$
	\begin{align}\label{kappa}
		|\kappa(p) - \kappa(\pi)| \leq \kappa_0 |p - \pi|.
	\end{align}
\end{assumption}
\begin{remark}\label{remark} Under Assumption~\ref{issODE}, for system $\dot{X}(t)=f(X(t),\kappa(X(t))+w(0,t))$, thanks to \cite[Lemma 4.6]{khalil2002nonlinear} we can prove the existence of constants $\sigma,M_{\sigma}, b_3 > 0$ such that the following inequality holds for $t \geq 0$
	\begin{align}
		|X(t)| \leq M_{\sigma}|X_0| e^{-\sigma t} + b_3 \operatorname{ess\ sup}_{0 \leq s \leq t} \|w(\cdot,s)\|_{\infty}.
	\end{align}
\end{remark}
\begin{remark}
	Under  Assumptions~\ref{globalLipschitz} and \ref{issODE}, using Gronwall's Lemma (see, e.g., Lemma A.1 in \cite{khalil2002nonlinear}), the following inequality holds
	\begin{align}
		M_4 \|(X, u)\| \leq \|(X, w)\| \leq M_3 \|(X, u)\|, \label{equivalence_constant}
	\end{align}	where
	\begin{align}
		\label{M3} M_3&=1+\kappa_0\max\{1,LD\}e^{LD},\\
		\label{M4}M_4&=\dfrac{1}{1+\kappa_0\max\{1,\kappa_0LD\}e^{LD(1+\kappa_0)}}.
	\end{align}
\end{remark}

	\subsection{Properties of the Quantizer}
The state $X$ of the plant and the actuator state $u$ are available only in quantized form. We consider here dynamic quantizers with an adjustable parameter of the form (see, e.g., \cite{bekiaris2020hybrid,brockett2000quantized,fkoudohode2024,liberzon2003hybrid,liberzon2006quantization})
\begin{equation}\label{quantizer}
	q_{\mu}(X,u)=( q_{1\mu}(X),q_{2\mu}(u))=\left(\mu q_1\left(\frac{X}{\mu}\right),\mu q_2\left(\frac{u}{\mu}\right)\right),
\end{equation} where $\mu>0$ can be manipulated and this is called zoom variable. The quantizers $q_1:\mathbb{R}^n\to \mathbb{R}^n$ and $q_2:\R\to \R$ are locally Lipchitz\footnote{ This assumption is required only for establishing well-posedness of the closed-loop system and it appears as it cannot be removed due to the potential non-measurability of the function $q_{2\mu}(u(y))$ in the integral in \eqref{predictor_quantizer} (see also, e.g., \cite{bekiaris2020hybrid,fkoudohode2024}).} functions that satisfy the following properties\\
P1: If $\|(X,u)\| \leq M$, then $\|(q_1(X)-X, q_2(u)-u)\| \leq \Delta$,\\
P2: If $\|(X,u)\|>M$, then $\|(q_1(X),q_2(u))\|>M-\Delta$,\\
P3: If $\|(X,u)\| \leq \hat{M}$, then $q_1(X)=0$ and $q_2(u)=0,$ \\
for some positive constants $M, \hat{M}$, and $\Delta$, with $M>\Delta$ and $\hat{M}<M$. 
\subsection{Quantized Predictor-Feedback Law}
The hybrid predictor-feedback law can be viewed as a quantized version of the predictor-feedback controller \eqref{nominalU}, in which the dynamic quantizer depends on a suitably chosen piecewise constant signal $\mu$. It is defined as 
\begin{equation}\label{control_quantizer}
U(t)=\left\{\begin{array}{ll}0, & 0 \leq t < t_{1}^* \\ \kappa( P_{\mu(t)}(t)), & t\ge t_{1}^*
\end{array}\right.,
\end{equation} with $P_{\mu}=p_{\mu}(D),$ where for $x\in[0,D]$ \begin{equation}\label{predictor_quantizer}
p_{\mu}(x)=q_{1\mu}(X)+\displaystyle\int_0^xf(p_{\mu}(y),q_{2\mu}(u(y)))dy.
\end{equation} 
The tunable parameter $\mu$ is selected as
\begin{equation}\label{switching_parameter}
\mu(t)= \begin{cases} 2\mathrm{e}^{2L(j+1) \tau} \mu_{0}& (j-1) \tau \leq t \leq j \tau+\bar{\tau}\delta_j, \\ & 1 \leq j \leq\left\lfloor\frac{t_{1}^*}{\tau}\right\rfloor,\\ \mu\left(t_{1}^*\right), & t \in\left[t_{1}^*, t_{1}^*+T\right.), \\ \Omega \mu\left(t_{1}^*+(i-2) T\right), & t \in\left[t_{1}^*+(i-1) T,\right. \\ & \left.t_{1}^*+i T\right), \quad i=2,3, \ldots\end{cases},
\end{equation} for some fixed, yet arbitrary, $\tau, \mu_0>0$, where $t_1^*=m\tau+\bar{\tau},$ for an $m\in\mathbb{Z}_+,$ $\bar{\tau}\in[0,\tau),$ and $\delta_m=1,\delta_j=0,j< m$, with $t_{1}^*$ being the first time instant at which the following holds 
\begin{align}
\nonumber&\left|\mu(t_{1}^*)q_1\left(\frac{X\left(t_{1}^*\right)}{\mu\left(t_{1}^*\right)}\right)\right|+\left\|\mu(t_{1}^*)q_2\left(\frac{u\left(\cdot, t_{1}^*\right)}{\mu\left(t_{1}^*\right)}\right)\right\|_{\infty}\\
&	\label{first_time_t0}\leq  (\overline{M}M- \Delta)\mu(t_1^*),
\end{align} where
\begin{align}
\label{M}
\overline{M} & = \frac{M_{4}}{M_3(1+M_0)}, \\
\label{M5}
M_5 & = \kappa_0 \max\{1, LD\} e^{LD},\\
\label{Omega}
\Omega & = \frac{M_5 \Delta (1+\lambda)(1+M_0)^2}{M_4 M}, \\
\label{T}
T & = -\frac{1}{\delta} \ln\left(\frac{\Omega}{1+M_0}\right),
\end{align} for some $\delta\in (0,\min\{\sigma,\nu\}),$
$\lambda$ is selected large enough in such a way that the following small-gain condition holds
\begin{equation}\label{small_gain}
\frac{b_3+1}{1+\lambda}<e^{-D},	
\end{equation} and $M_0$ is defined such that
\begin{align}
\nonumber& M_{0}=\left(1-\phi\right)^{-1} \left(1-\varphi_1\right)^{-1}\max \left\{e^{D(\nu+1)} ;\phi M_{\sigma}\right\}\\
&+(1-\varphi_1)^{-1}\max \left\{M_{\sigma}; (1+\varepsilon) \left(1-\phi\right)^{-1}e^{D(\nu+1)}b_3\right\},
\end{align} where $0<\phi<1$ and $0<\varphi_{1}<1$ with
\begin{align} \label{phi and psy}
\phi= \frac{1+\varepsilon}{1+\lambda}e^{D\left(\nu+1\right)} \text{ and }\varphi_{1}=(1+\varepsilon)(1-\phi)^{-1}\phi b_3,
\end{align} for some  $\varepsilon>0$. The choice of $\nu,\varepsilon$ guarantees that $\phi<1,\varphi_1<1$, which is always possible given \eqref{small_gain} (see also the proof of Lemma~\ref{Lemma2} in Section~\ref{stabStateQuantization}). 
\section{Stability Under State Quantization}\label{stabStateQuantization}
\begin{theorem}\label{Theorem1}
	Consider the closed-loop system consisting of the plant \eqref{pde_representation}--\eqref{pde_representation1} and the switched predictor-feedback law \eqref{control_quantizer}--\eqref{switching_parameter}. Under Assumptions~\ref{globalLipschitz} and \ref{issODE}, if $\Delta$ and $M$ satisfy  \begin{equation}\label{conditionMDelta}
		\frac{\Delta}{M}<\frac{M_4}{(1+M_0)\max \{M_5(1+\lambda)(1+M_0),2M_5\}},
	\end{equation}	then for all $X_{0} \in \mathbb{R}^{n}$, $u_{0} \in \mathcal{C}_{\rm rpw}([0, D], \mathbb{R})$, there exists a unique solution such that $X(t) \in AC\left(\mathbb{R}_{+}, \mathbb{R}^{n}\right)$, for each $t \in \mathbb{R}_{+}$ $u(\cdot, t) \in \mathcal{C}_{\rm rpw}([0,D], \mathbb{R})$, and for each $x \in[0,D]$ $u(x,\cdot) \in \mathcal{C}_{\rm rpw}\left(\mathbb{R}_{+}, \mathbb{R}\right)$, which satisfies
	\begin{align}
		&|X(t)|+\left\|u(\cdot,t)\right\|_{\infty}
		\label{stability_result_u} \leq  \gamma\left(|X_0|+\left\|u_{0}\right\|_{\infty}\right)^{\left(2-\frac{\ln \Omega}{T} \frac{1}{L}\right)} \mathrm{e}^{\frac{\ln \Omega}{T}t},
	\end{align}	where
\begin{align}
	\nonumber	\gamma&=\frac{2}{M_4}\max \left\{\frac{M_4M\mu_{0}}{\Omega} e^{2L \tau} , M_3\right\}  \max \left\{\frac{1}{\mu_0(M \overline{M}-2 \Delta)}, \right. \nonumber \\
	&\quad\quad \left. 1 \right\} \left(\frac{1}{\mu_0(M \overline{M}-2 \Delta)}\right)^{\left(1-\frac{\ln \Omega}{T} \frac{1}{L}\right)}. \label{gamma}
\end{align}
\end{theorem} 
The proof relies on the following two lemmas.
\begin{lemma}\label{Lemma1}
	Let $\Delta$ and $M$ satisfy \eqref{conditionMDelta}, 
	there exists a time $t_{1}^*$ satisfying 
	\begin{equation}\label{t0}
		t_{1}^* \leqslant \frac{1}{L} \ln\left(\dfrac{\left|X_{0}\right|+\left\|u_{0}\right\|_{\infty}}{\mu_{0}(M \overline{M}-2 \Delta)}\right),
	\end{equation} 
	such that \eqref{first_time_t0} holds, and thus, the following also holds 
	\begin{equation}\label{bound_in_time_t0}
		|X(t_{1}^*)|+\|u(\cdot, t_{1}^*)\|_{\infty} \leq  \overline{M}M\mu(t_{1}^*).
	\end{equation}
\end{lemma}
\begin{proof}	
	For $0 \leq t < t_1^*$, with \eqref{control_quantizer} one has $U(t)=0$. Thus, the open-loop system is given by
	\begin{align}
		\dot{X}(t) &= f(X(t), u(0, t)), \label{pde_representation0} \\
		u_t(x, t) &= u_x(x, t), \label{pde_representation4} \\
		u(D, t) &= 0. \label{pde_representation3}
	\end{align}
	Using the method of characteristics, the solution to the transport equation is $u(x, t) = u_0(x+t)$ for $0 \leq x+t \leq D$ and $u(x, t) = 0$ for $x+t > D$. Hence, 
	\begin{equation}\label{estimation_u_norm}
		\|u(\cdot, t)\|_{\infty} \leq \|u_0\|_{\infty}.
	\end{equation}
	By Assumption~\ref{globalLipschitz} and Gronwall's Lemma, one has
	\begin{equation}\label{estimation_X_norm}
		|X(t)| \leq e^{Lt} \left(|X_0| + \|u_0\|_{\infty}\right).
	\end{equation}
	Combining \eqref{estimation_u_norm} and \eqref{estimation_X_norm}, we obtain
	\begin{align}
		|X(t)| + \|u(\cdot, t)\|_{\infty} &\leq 2e^{Lt}\left(|X_0| +  \|u_0\|_{\infty}\right).\label{normXu}
	\end{align}
	The rest of the proof follows in exactly the same way as \cite[Lemma 1]{fkoudohode2024}.
\end{proof}

\begin{lemma}\label{Lemma2}
	Select $\lambda$ large enough in such a way that the small-gain condition \eqref{small_gain} holds. Then the solution to the target system \eqref{targetsystem_without_quantizer}--\eqref{targetsystem_without_quantizer3} 	with the quantized controller \eqref{control_quantizer}, 
		resulting in $w(D,t)=\kappa\left(P_{\mu}(t)\right))-\kappa\left(P(t)\right)$ with $u$ given in terms of $(X,w)$ by the inverse backstepping transformation \eqref{backstepping_direct_transformation}, \eqref{pixt}, which verify, for fixed $\mu$,
	\begin{equation}\label{hyplemma2}
		|X(t_{1}^*)|+\|w( \cdot, t_{1}^*)\|_{\infty}\le M_{3}\overline{M}M\mu, 
	\end{equation} they satisfy for $t_{1}^*\le t< t_{1}^*+T$
	\begin{align}
		\nonumber& |X(t)|+\|w(\cdot,t)\|_{\infty}\leqslant \max\left\{ M_{0} e^{-\delta (t-t_{1}^*)}\left(\left|X(t_{1}^*)\right| \right.\right.\\
		&\left. +\left\|w( \cdot, t_{1}^*)\right\|_{\infty}\right),\left. \Omega M_{3}\overline{M}M\mu  \right\}.\label{norm_X_w} 
	\end{align} In particular, the following holds
	\begin{align}\label{normXu1}
		|X(t_{1}^*+T)|+\|w(\cdot,t_{1}^*+T)\|_{\infty} \leq \Omega M_{3}\overline{M}M\mu. 
	\end{align}
\end{lemma}\begin{proof}				
	From \eqref{small_gain}, note that the function
	\begin{equation}
		h\left(s_1, s_2\right) = \frac{1+s_1}{1+\lambda} e^{D\left(s_2+1\right)}\left(b_3(s_1+1)+1\right),
	\end{equation}
	is continuous at $(0, 0)$ and satisfies $h(0,0) < 1$. Consequently, there exist constants $\varepsilon > 0$ and $\nu > 0$ such that $h\left(\varepsilon, \nu\right) < 1$, that is,
	\begin{equation}\label{small_gain2}
		\frac{1+\varepsilon}{1+\lambda} e^{D(\nu+1)}\left(b_3(\varepsilon+1)+1\right) < 1.
	\end{equation}
	Furthermore, this condition implies 
	\begin{equation}\label{small_gain1}
		\frac{1+\varepsilon}{1+\lambda} e^{D\left(\nu+1\right)} < 1.
	\end{equation}
	Using the fading memory lemma \cite[Lemma 7.1]{karafyllis2021input}, it follows from Remark~\ref{remark} that for any $\varepsilon > 0$, there exists $\sigma_1 \in (0,\sigma)$ such that:
	\begin{align}
		\nonumber	&|X(t)|e^{\sigma_1(t-t^*)} \leqslant M_{\sigma}\left|X(t_1^*)\right|\\
		&+(1+\varepsilon)b_3 \operatorname{ess\ sup}_{t_{1}^* \leq s \leq t}\left(\|w(\cdot,s)\|_{\infty} e^{\sigma_1(s-t^*)}\right). \label{eqXw}
	\end{align}
	
	Next, consider the transport subsystem \eqref{targetsystem_without_quantizer2}, \eqref{targetsystem_without_quantizer3}. Applying the ISS estimate in the sup-norm from \cite[estimate (2.23)]{KKnonlocal} (see also \cite[estimate (3.2.11)]{karafyllis2021input}), for any $\nu > 0$, we have
	\begin{align}
		\nonumber	\|w(\cdot,t)\|_{\infty} &\leq e^{-\nu(t-t_{1}^*-D)}e^{D}\left\|w(\cdot,t_{1}^*)\right\|_{\infty}\\
		\label{wmu1d}	&+e^{D\left(1+\nu\right)} \operatorname{ess\ sup}_{t_{1}^* \leq s \leq t}(|d(s)|),
	\end{align} where 
	\begin{align}
		\label{d1}	d(t)&=\kappa\left(p_{\mu}(D,t)\right) -\kappa\left(p(D,t)\right).
	\end{align} By the fading memory inequality \cite[Lemma 7.1]{karafyllis2021input}, there exists $\delta_2 \in (0, \nu)$ such that
	\begin{align}
		\nonumber	\|w(\cdot,t)\|_{\infty} e^{\delta_2 (t-t_{1}^*)} &\leqslant e^{D\left(\nu+1\right)}\left\|w(\cdot,t_{1}^*)\right\|_{\infty}+e^{D\left(\nu+1\right)}(1+\varepsilon) \\&\times\operatorname{ess\ sup}_{t_{1}^* \leq s \leq t}\left(|d(s)| e^{\delta_2 (s-t_{1}^*)}\right). \label{eqWd}
	\end{align}
	
	We use next equations \eqref{pxt}, \eqref{predictor_quantizer}, and the global Lipschitzness assumptions (Assumptions~\ref{globalLipschitz} and \ref{issODE}) on $\kappa$ and $f$ to obtain
	\begin{align}
		\nonumber &|p(x)-p_{\mu}(x)|
	\le \left|q_{1\mu}(X)-X\right| +L\displaystyle\int_{0}^{x} \left|p(\xi)-p_{\mu}(\xi)\right| d\xi\\
		\nonumber	&+L\displaystyle\int_{0}^{x} |q_{2\mu}(u(\xi))-u(\xi)| d\xi\\
		\nonumber	&\le \left|q_{1\mu}(X)-X\right| +LD\|q_{2\mu}(u)-u\|_{\infty}\\
		&+L\displaystyle\int_{0}^{x} \left|p(\xi)-p_{\mu}(\xi)\right| d\xi, 
	\end{align} and hence, by \cite[Lemma A.1]{khalil2002nonlinear} we get 
	\begin{align}
	\label{Pmutpt}	 |p(x)-p_{\mu}(x)| &\le \max\{1,LD\}e^{LD}\left(\left|q_{1\mu}(X)-X\right|\right.\\
	\nonumber	 &\left.+\|q_{2\mu}(u)-u\|_{\infty}\right). 
	\end{align}
	Therefore, from \eqref{kappa}, \eqref{d1}, and \eqref{Pmutpt} we obtain
		\begin{align}
		|d|&\le M_5\left(\left|q_{1\mu}(X)-X\right| +\|q_{2\mu}(u)-u\|_{\infty}\right),
	\end{align} where $M_5$ is defined in \eqref{M5} and $u$ is given in terms of $(X,w)$ by the inverse backstepping transformation \eqref{backstepping_direct_transformation}, \eqref{pixt}. Provided that 
	\begin{align}
		\label{condition}
		\frac{	\Omega M_{4}M\mu}{(1+M_0)^2}  \leq|X|+\|w\|_{\infty} \leqslant M_{4} M\mu,
	\end{align}
	thanks to the property $\rm P1$ of the quantizer, the left-hand side of bound \eqref{equivalence_constant}, and the definition \eqref{Omega}, we obtain 
	\begin{align}
		\nonumber	|d|& \le  M_5\Delta\mu\\&\leqslant \frac{(1+M_0)^2 M_5\Delta}{\Omega M_4M} \left(|X|+\|w\|_{\infty}\right) \\
		& \leqslant \frac{1}{1+\lambda}\left(|X|+\|w\|_{\infty}\right). 
	\end{align}	 
	Let us define $\delta$ such that $\delta \in (0, \min\{\sigma_1,\delta_2\})$. Now, for $t_{1}^*\le t <t_1^*+T$, let us define the following quantities   
	\begin{align}
		\label{defnormw}\|w\|_{[t_{1}^*, t]}&:=\operatorname{ess\ sup}_{t_{1}^* \leqslant s \leqslant t}\|w(\cdot, s)\|_{\infty}e^{\delta (s-t_{1}^*)} ,\\
		|\label{defnormX}X|_{[t_{1}^*, t]}&:=\operatorname{ess\ sup}_{t_{1}^* \leqslant s \leqslant t}|X(s)| e^{\delta (s-t_{1}^*)}.
	\end{align}
	
	As long as the solutions satisfy \eqref{condition} we get
	\begin{align}\label{estimate_d}
		\|d\|_{[t_{1}^*, t]}	&\leqslant \frac{1}{1+\lambda}\|w\|_{[t_{1}^*, t]}+\frac{1}{1+\lambda}|X|_{[t_{1}^*, t]}.
	\end{align}
	From now on the proof is identical to \cite[Lemma 2]{fkoudohode2024}.
\end{proof}
{\em Proof of Theorem~\ref{Theorem1}:} We proceed in the same manner as in \cite[Proof of Theorem~1]{fkoudohode2024} applying Lemma~\ref{Lemma1} and Lemma~\ref{Lemma2}.

We now prove the well-posedness of the system. In the interval $\left[0, t_{1}^*\right)$, where no control is applied, the system is described by \eqref{pde_representation0}--\eqref{pde_representation3}. The existence and uniqueness of solutions within this interval are guaranteed by the global Lipschitzness assumption on $f$ and the explicit solution of the transport subsystem \eqref{pde_representation3} via the method of characteristics. Moreover, we have $X(t) \in AC\left([0, t_{1}^*), \mathbb{R}^{n}\right)$ and $u \in \mathcal{C}_{\rm rpw}([0, D] \times [0, t_{1}^*), \mathbb{R})$. For $t > t_{1}^*$, the system described by \eqref{pde_representation}--\eqref{pde_representation1}, along with the quantized controller $U$ defined in \eqref{control_quantizer}, satisfies the conditions outlined in \cite[Theorem 8.1]{karafyllis2021input}. Specifically, the terms $F(X, u) = f(X, u)$ and $\varphi(\mu, u, X) = U(\mu, u, X)$ fulfill these assumptions. In particular, the control mapping $U(\mu,\cdot,\cdot):L^{\infty}([0,D];\R)\times \R^n\to \R$ defined in \eqref{control_quantizer}, \eqref{predictor_quantizer} is locally Lipschitz on bounded sets, owing to the Lipschitz continuity of $q_1$ and $q_2,$ and the Lipschitzness assumption (Assumptions~\ref{globalLipschitz} and \ref{issODE}) on $f$ and $\kappa$. The rest of the proof follows in the exact same manner as the respective proof in \cite{fkoudohode2024}. 
\section{Extension to Input Quantization}\label{inputquantization}
In this section, we address the case of input quantization. Here, measurements of the states are assumed to be available, and the control input is defined as follows
\begin{equation}\label{control_quantizerinput}
	U(t) = 
	\begin{cases} 
		0, & 0 \leq t < \bar{t}_{1}^*, \\  
		\bar{q}_{\mu}\big(U_{\rm nom}(t)\big), & t \ge \bar{t}_{1}^*,
	\end{cases}
\end{equation}
where \( U_{\rm nom}(t) \) is specified in \eqref{nominalU}. The quantizer is a locally Lipschitz function \( \bar{q}_{\mu}: \mathbb{R} \to \mathbb{R} \), defined as \( \bar{q}_{\mu}(\bar{U}) = \mu \bar{q}\big(\frac{\bar{U}}{\mu}\big) \), and it satisfies the following properties
\begin{itemize}
	\item[\(\bar{\rm P}1\):] If \( |\bar{U}| \leq M \), then \( |\bar{q}(\bar{U}) - \bar{U}| \leq \Delta \),
	\item[\(\bar{\rm P}2\):] If \( |\bar{U}| > M \), then \( |\bar{q}(\bar{U})| > M - \Delta \),
	\item[\(\bar{\rm P}3\):] If \( |\bar{U}| \leq \hat{M} \), then \( \bar{q}(\bar{U}) = 0 \).
\end{itemize}			
The tunable parameter $\mu$ is selected as 
\begin{equation}\label{switching_parameterinput}
	\mu(t)= \begin{cases} 2\mathrm{e}^{2L (j+1)\tau} \mu_{0}, & (j-1) \tau \leq t \leq j \tau+\bar{\tau}_1\delta_j, \\ & 1 \leq j \leq\left\lfloor\frac{\bar{t}_{1}^*}{\tau}\right\rfloor,\\ \mu\left(\bar{t}_{1}^*\right), & t \in\left[\bar{t}_{1}^*, \bar{t}_{1}^*+T\right), \\ \Omega \mu\left(\bar{t}_{1}^*+(i-2) T\right), & t \in\left[\bar{t}_{1}^*+(i-1) T,\right. \\ & \left.\bar{t}_{1}^*+i T\right), \quad i=2,3, \ldots\end{cases},
\end{equation} for some fixed, yet arbitrary, $\tau, \mu_0>0$, where $\bar{t}_1^*=\overline{M}\tau+\bar{\tau}_1,$ for an $\overline{M}\in\mathbb{Z}_+,$ $\bar{\tau}_1\in[0,\tau),$ and $\delta_{\overline{M}}=1,\delta_j=0,j< \overline{M}$, with $\bar{t}_{1}^*$ being the first time instant at which the following holds  
\begin{equation}\label{bound_in_time_t0input}
	|X\left(\bar{t}_{1}^*\right)| + \|u\left(\cdot,\bar{t}_1^*\right)\|_{\infty} \leq \frac{M \overline{M}}{M_5} \mu(\bar{t}_{1}^*).
\end{equation}
The parameters in \eqref{switching_parameterinput} and \eqref{bound_in_time_t0input} are defined in \eqref{M3}--\eqref{T}.
\begin{theorem}\label{Theorem2}
	Consider the closed-loop system comprising the plant \eqref{pde_representation}--\eqref{pde_representation1} and the switched predictor-feedback law \eqref{control_quantizerinput}, \eqref{switching_parameterinput} with \eqref{nominalU}. Under Assumptions~\ref{globalLipschitz} and \ref{issODE}, if \( \Delta \) and \( M \) satisfy $\frac{\Delta}{M} < \frac{M_4}{M_5 (1 + \lambda)(1 + M_0)^2},$
	then for all \( X_{0} \in \mathbb{R}^{n} \), \( u_{0} \in \mathcal{C}_{\rm rpw}([0, D], \mathbb{R}) \), there exists a unique solution such that \( X(t) \in AC(\mathbb{R}_{+}, \mathbb{R}^{n}) \), \( u(\cdot, t) \in \mathcal{C}_{\rm rpw}([0,D], \mathbb{R}) \) for each \( t \in \mathbb{R}_{+} \), \( u(x,\cdot) \in \mathcal{C}_{\rm rpw}(\mathbb{R}_{+}, \mathbb{R}) \) for each \( x \in [0,D] \),
	which satisfies
	\begin{equation}\label{stability_result_uinput}
		|X(t)| + \|u(\cdot,t)\|_{\infty} \leq \bar{\gamma}\big(|X_0| + \|u_{0}\|_{\infty}\big)^{2 - \frac{\ln \Omega}{T} \frac{1}{L}} \mathrm{e}^{\frac{\ln \Omega}{T}t},
	\end{equation}
	where
	\begin{align}
		\bar{\gamma} &= \frac{2}{M_4} \max \bigg\{\frac{M_4 M}{\Omega M_5} \mathrm{e}^{2L \tau} \mu_{0}, M_3\bigg\} \max \bigg\{\frac{M_5}{\mu_0 M \overline{M}}, 1\bigg\} \nonumber \\
		&\quad \times \bigg(\frac{M_5}{\mu_0 M \overline{M}}\bigg)^{1 - \frac{\ln \Omega}{T} \frac{1}{L}}.
	\end{align}
\end{theorem}
	The proof of Theorem~\ref{Theorem2} relies on the following two lemmas, which are analogous to the case of state quantization.
\begin{lemma}\label{Lemma3}
	For the time interval \( [0, \bar{t}_{1}^*) \) where \( U(t) = 0 \), the following inequality holds
	\begin{equation}
		|X(t)| + \|u(\cdot,t)\|_{\infty} \leq 2\mathrm{e}^{L t}\big(|X_{0}| + \|u_{0}\|_{\infty}\big).
	\end{equation}
\end{lemma}
\begin{proof}
	Identical to the proof of Lemma~\ref{Lemma1}.
\end{proof}
\begin{lemma}\label{Lemma4}
	Select $\lambda$ large enough in such a way that the small-gain condition \eqref{small_gain} holds. Then the solutions to the target system \eqref{targetsystem_without_quantizer}--\eqref{targetsystem_without_quantizer3} with the quantized controller \eqref{nominalU}, \eqref{control_quantizerinput}, \eqref{switching_parameterinput}, which verify, for fixed $\mu$,
	\begin{equation}\label{hyplemma4input}
		|X(\bar{t}_{1}^*)|+\|w(\cdot,\bar{t}_{1}^*)\|_{\infty}\le \frac{M_4M\mu}{(1+M_0)M_5},
	\end{equation} they satisfy for $\bar{t}_{1}^*\le t< \bar{t}_{1}^*+T$
	\begin{align}
		\nonumber& |X(t)|+\|w(\cdot,t)\|_{\infty}\leqslant \max\left\{ M_{0} e^{-\delta (t-\bar{t}_{1}^*)}\left(\left|X(\bar{t}_{1}^*)\right| \right.\right.\\
		&\left. +\left\|w(\cdot,\bar{t}_{1}^*)\right\|_{\infty}\right),\left.  \frac{\Omega M_4M\mu }{(1+M_0)M_5}  \right\}.\label{norm_X_winput}
	\end{align} In particular, the following holds
	\begin{align}\label{normXu1input}
		|X(\bar{t}_{1}^*+T)|+\|w(\cdot,\bar{t}_{1}^*+T)\|_{\infty} \leq  \frac{\Omega M_4M\mu}{(1+M_0)M_5}.
	\end{align}
\end{lemma}
\begin{proof}
	For $\bar{t}_{1}^*\le t< \bar{t}_{1}^*+T,$ the system is defined by \eqref{pde_representation}--\eqref{pde_representation1} under the switched predictor-feedback law \eqref{nominalU}, \eqref{control_quantizerinput}, \eqref{switching_parameterinput}. Using the same strategy, as in the proof of  Lemma~\ref{Lemma2}, i.e., combining Remark~\ref{remark}, the ISS estimate in sup-norm in \cite[estimate (2.23)]{KKnonlocal}, and the fading memory lemma  \cite[Lemma 7.1]{karafyllis2021input}, for every $\nu, \varepsilon>0$ satisfying \eqref{small_gain2}, there exists $\delta \in (0, \min\{\sigma,\nu\})$ such that, using the definitions \eqref{defnormw} and \eqref{defnormX}, the following inequalities hold
	\begin{align}\label{norm_x_winput}
		|X|_{[\bar{t}_{1}^*, t]} \leqslant M_{\sigma}|X(\bar{t}_{1}^*)|+(1+\varepsilon) b_3\|w\|_{[\bar{t}_{1}^*, t]},
	\end{align}	and
	\begin{align}					\nonumber\|w\|_{[\bar{t}_{1}^*, t]}&\le e^{D\left(\nu+1\right)}(1+\varepsilon) \operatorname{ess\ sup}_{\bar{t}_{1}^*\leqslant s \leqslant t}\left(\left|\bar{d}(s)\right| e^{\delta (s-\bar{t}_{1}^*)}\right) \label{norm_w_0_tinput} \\
		&+e^{D}\left\|w(\cdot,\bar{t}_{1}^*)\right\|_{\infty},
	\end{align}	where
	\begin{align}
		\bar{d}(t)=U_{\rm nom}(t)-\mu(t)\bar{q}\left(\dfrac{U_{\rm nom}(t)}{\mu(t)}\right) ,	\label{d1input}
	\end{align} with $ U_{\rm nom}$ and $\mu$ given in \eqref{nominalU} and \eqref{switching_parameterinput}, respectively. Next, let us proceed to approximate for fixed $\mu$ the term\\ $\displaystyle\operatorname{ess\ sup}_{\bar{t}_{1}^* \leqslant s\leqslant t}\left(\left|\bar{d}(s)\right| e^{\delta (s-\bar{t}_{1}^*)}\right)$ in \eqref{norm_w_0_tinput}.  
	Provided that 
	\begin{align}
		\label{conditioninput}
		\frac{\Omega MM_4\mu}{(1+M_0)^2M_5} \leq|X|+\|w\|_{\infty} \leqslant \dfrac{M_{4} M\mu}{M_5},
	\end{align} 
	using the property $\rm \bar{P}1$ of the quantizer, the left-hand side of bound \eqref{equivalence_constant}, and the definition \eqref{Omega}, we obtain \begin{align}
		\nonumber	\left|\bar{d}\right|&=\mu\left|\bar{q}\left(\dfrac{U_{\rm nom}(t)}{\mu}\right)-\dfrac{U_{\rm nom}(t)}{\mu}\right|   \\
		& \leqslant \frac{1}{1+\lambda}\left(|X|+\|w\|_{\infty}\right), 
	\end{align}	where we also used the fact that $|U_{\rm nom}|\le M_5(|X|+\|u\|_{\infty})$ which follows from \eqref{kappa} and the fact that $\|p\|_{\infty}\le \max\{1,LD\}e^{LD}\left(|X|+\|u\|_{\infty}\right),$ that is obtained using \eqref{pxt}, \eqref{lipschitzf}, and Gronwall's Lemma.
	Hence, using \eqref{norm_w_0_tinput} and the fact that \eqref{small_gain1} holds we obtain
	\begin{align}
		\|w\|_{[\bar{t}_{1}^*, t]}\leqslant\left(1-\phi\right)^{-1}e^D\left\|w(\cdot,\bar{t}_{1}^*)\right\|_{\infty}+\left(1-\phi\right)^{-1}\phi |X|_{[\bar{t}_{1}^*, t]} ,\label{normwinput}
	\end{align}
	with $\phi=\frac{1+\varepsilon}{1+\lambda}e^{D\left(\nu+1\right)}<1.$
	Combining the inequalities \eqref{norm_x_winput} and \eqref{normwinput}, thanks to the definitions \eqref{defnormw}, \eqref{defnormX}, and to the small-gain condition \eqref{small_gain}, repeating the respective arguments from the proof of Lemma~\ref{Lemma2}, we arrive at
	\begin{equation}\label{norm_X_u1input}
		|X(t)|+\|w(\cdot,t)\|_{\infty}\leqslant M_0e^{-\delta (t-\bar{t}_{1}^*)}\left(|X(\bar{t}_{1}^*)|+\left\|w(\cdot,\bar{t}_{1}^*)\right\|_{\infty}\right). 
	\end{equation}
	The rest of the proof utilizes the same reasoning as in \cite[Lemma 4]{fkoudohode2024}.
\end{proof}	
{\em Proof of Theorem~\ref{Theorem2}:} Identical to the proof of Theorem~\ref{Theorem1}.
\section{Conclusions and current work}\label{conclusion}
We have established global asymptotic stability for nonlinear systems with input delay, addressing both state and input quantization, through the introduction of a nonlinear, switched predictor-feedback control law. Our proof utilized a combination of the backstepping method, small-gain techniques, and input-to-state stability arguments. Our current efforts focus on relaxing the global Lipschitzness assumptions on the system dynamics and feedback law achieving semi-global stabilization results, as well as on relaxing the exponential stability assumption of the nominal, delay/quantization-free closed-loop system.

\bibliography{ifacconf} 

\end{document}